\documentclass[12pt]{article}

\usepackage{alltt, amsmath, amssymb,amsthm, amsfonts, graphicx, verbatim}

\newtheorem{theorem}{Theorem}
\newtheorem{lemma}[theorem]{Lemma}
\newtheorem{corollary}[theorem]{Corollary}
\newtheorem{proposition}[theorem]{Proposition}

\newtheorem{assumption}{Assumption}
\newtheorem{remark}[theorem]{Remark}

\newcommand{\bea}{\begin{eqnarray}}
\newcommand{\be}{\begin{equation}}
\newcommand{\ba}{\begin{array}}
\newcommand{\bl}{\begin{lemma}}
\newcommand{\bc}{\begin{corollary}}
\newcommand{\bt}{\begin{theorem}}
\newcommand{\br}{\begin{remark}}
\newcommand{\bp}{\begin{proof}}
\newcommand{\bi}{\begin{itemize}}
\newcommand{\bprop}{\begin{proposition}}

\newcommand{\eea}{\end{eqnarray}}
\newcommand{\ee}{\end{equation}}
\newcommand{\ea}{\end{array}}
\newcommand{\el}{\end{lemma}}
\newcommand{\ec}{\end{corollary}}
\newcommand{\et}{\end{theorem}}
\newcommand{\er}{\end{remark}}
\newcommand{\ep}{\end{proof}}
\newcommand{\ei}{\end{itemize}}
\newcommand{\eprop}{\end{proposition}}

\newcommand{\nn}{\nonumber}
\newcommand{\ra}{\rightarrow}
\newcommand{\lra}{\longrightarrow}

\title{A vector minmax problem for controlled Markov chains}

\date{}

\begin{document}

\author{Sameer Kamal\footnote{School of Technology and Computer Science,
Tata Institute of Fundamental Research, Homi Bhabha Road,
Mumbai-400005, India. E-mail: sameer.kamal@gmail.com. This work was
supported in part by an Infosys Fellowship.}}

\maketitle

\begin{abstract}
The problem of controlling a finite state Markov chain in the
presence of an adversary so as to ensure desired performance levels
for a vector of objectives is cast in the framework of Blackwell
approachability. Relying on an elementary two time scale
construction a control scheme is proposed which ensures almost sure
convergence to the desired set regardless of the adversarial
actions.
\end{abstract}

\noindent \textbf{Key words:} controlled Markov chains, Blackwell
approachability, two time scales, stationary strategies,
multi--objective optimization

\section{Introduction}

Many control problems in practice have two features that put them
outside of the classical framework of deterministic or stochastic
optimal control theory: presence of unknown disturbances and
multiple objectives. One common approach for addressing the former
issue is to treat the disturbances as actions of an adversary and
plan against the worst case scenario thereof. This makes the problem
a two person zero sum game. While the classical two person zero sum
stochastic games are fully analyzable through the associated Shapley
equation, this is not the case when there are many objectives. In a
seminal article, Blackwell~\cite{blackwell} provided a framework for
addressing this `vector minmax' problem in case of repeated games,
providing both the necessary and sufficient conditions for
attainability of the objectives (what came to be known as
\textit{Blackwell approachability}) and a scheme for achieving the
same. This is becoming a popular model for addressing engineering
problems with aforementioned features, see, e.g., Hou et
al~\cite{hou} for a recent application. The framework has also found
application in strategic learning literature in economics and
computer science, see, e.g., Young~\cite{young}. As observed above,
many engineering situations call for going beyond the repeated game
model and consider a controlled Markov dynamics instead. In an
important work, Shimkin and Shwartz~\cite{shimkin} studied this
problem for controlled Markov chains and proposed a scheme to ensure
Blackwell approachability. Their scheme depends on updating
strategies at return times to a fixed state, which allows them to
exploit the regenerative nature of such visits. This is necessitated
by the fact that there appears to be a need to hold the policy fixed
for some time -- the interval between two return times in their case
-- for the `learning' to take place. For a large chain, the return
times can be infrequent, rendering the convergence slower. Motivated
by this, we propose an alternative scheme here that holds a policy
constant for durations that are short initially and can become
longer gradually, thus capturing the `exploration-exploitation'
trade--off. Each choice of strategy is associated with a positive
re--scaled time duration and whenever the player switches to a new
strategy he retains it for the associated re--scaled duration of
time. Almost sure convergence of the running average cost to the
desired set is then established under standard conditions. A major
ingredient in our proof is an elementary two time scale argument and
the proposed scheme is designed to exploit the two time scale
feature in an essential way.

The paper is organized as follows. Section~\ref{basic} describes the
problem set--up and introduces the notation and some preliminary
concepts. Section~\ref{tts} develops an elementary two time scale
result which plays a crucial role in the proof of convergence and
around which our scheme is built in the first place.
Section~\ref{asc} proves the main convergence result,
Theorem~\ref{theoremmain}. Section~\ref{conclusion} concludes by
outlining some further possibilities.

\section{Basic setup}

\label{basic}

\emph{The model.}\ \  Consider a system evolving as a controlled
Markov chain on a finite state space $S$ with a reward associated
with each transition. We assume that the reward is always some
vector from a compact set $K\subset \mathbb R^d$. Let $U^p$ and
$U^a$ be finite action spaces. Let $(\theta_n)$ denote the
aforementioned controlled Markov chain on $S$ with transition kernel
$p(\theta^\prime| \theta, u^p, u^a)$ for $\theta^\prime, \theta \in
S, u^p \in U^p, u^a \in U^a$. Let $\mathcal P(U^p)$ denote the set
of probability distributions on the space $U^p$. Let $\Pi^p$ denote
the set of all maps, or strategies, from $S$ to $\mathcal P(U^p)$.
Similarly, let $\Pi^a$ denote the set of all strategies from $S$ to
$\mathcal P(U^a)$. Depending on the past the player and the
adversary \textit{independently} choose their current strategies
from $\Pi^p$ and $\Pi^a$ respectively. Let $(u^p_n)$, $(u^a_n)$ be
the actual control sequences chosen by the player and the adversary
from $U^p$, $U^a$ respectively. At time step $n$ the one step reward
is given by $\kappa\left(\theta_n, u^p_n, u^a_n\right)$. Let $x_n$
denote the vector for current average reward. The iterative equation
for the average reward becomes
\[
x_{n+1} = x_{n} + 1/(n+1)[\kappa\left(\theta_n, u^p_n, u^a_n\right)
- x_n].
\]

\emph{Main goal.}\ \  The aim of the main player is to have the
average reward asymptotically approach a certain desirable subset
$D\subset K \left(\subset \mathbb R^d\right)$ by suitably choosing
his strategy at each step. More precisely, the player seeks to
choose his sequence of strategies in such a manner that no matter
what sequence of strategies the adversary chooses, with probability
one all limit points of the sequence $(x_n)$ lie in $\bar D$ where
$\bar D$ denotes the closure of $D$.

\emph{Assumptions.}\ \  In our analysis we restrict our attention to
the case where $\bar D$ is convex. However, see
Section~\ref{conclusion} for possible extension to the case of
non--convex $\bar D$. Next,  assume that when the strategies for the
main player and the adversary are held fixed at arbitrary strategies
$\pi^p\in\Pi^p$ and $\pi^a\in\Pi^a$ respectively then the Markov
chain $(\theta_n)$ is ergodic. Let $\eta^{(\pi^p,\pi^a)}(\cdot)$
denote the corresponding stationary measure on state space $S$ with
the strategies for the player and the adversary held fixed. Define
the corresponding average reward $\bar \kappa(\pi^p,\pi^a)$ as
\[
\bar \kappa(\pi^p,\pi^a):=\sum_{\theta \in S}\sum_{u^p\in
U^p}\sum_{u^a\in U^a}\kappa(\theta, u^p,
u^a)\eta^{(\pi^p,\pi^a)}(\theta)\pi^p(u^p|\theta)\pi^a(u^a|\theta).
\]

For any point $x$, let $x_{\bar D}$ be the (unique) point in $\bar
D$ closest to $x$. For the rest of this paper we work under the
following assumption which is standard for Blackwell
approachability:

\begin{assumption}
\label{mainassumption} For every $x\in K\backslash \bar D$ there
exists a player strategy $\pi_x^p$ satisfying the following
inequality:
\[
\inf_{\pi^a\in\Pi^a}\langle \bar \kappa(\pi^p_x,\pi^a)-x_{\bar D},
x_{\bar D} -x\rangle> 0.
\]
In words, the hyperplane through $x_{\bar D}$ perpendicular to the
line segment $xx_{\bar D}$ separates $x$ from the set $\{\bar
\kappa(\pi^p_x,\pi^a):\pi^a\in\Pi^a\}$.
\end{assumption}

For $\rho\in\mathbb R_+$, let $B(x, \rho)$ denote the open ball of
radius $\rho$ centered at $x$.

\bl There exists a map $\rho(\cdot) : K\backslash \bar D \lra
\mathbb R_+$, such that for any $x\in K\backslash \bar D$, we have
\be \label{satisfactory} \inf_{y\in B(x,
\rho(x))}\inf_{\pi^a\in\Pi^a}\langle \bar
\kappa(\pi^p_x,\pi^a)-y_{\bar D},y_{\bar D}-y\rangle > 0. \ee \el

\bp

For any $x\in K\backslash \bar D$, by
Assumption~\ref{mainassumption} there exists a player strategy
$\pi^p_x$ and an $\epsilon > 0$ such that
\[
\inf_{\pi^a\in\Pi^a}\langle \bar \kappa(\pi^p_x,\pi^a)-x_{\bar D},
x_{\bar D} -x\rangle> \epsilon.
\]
 We get

\bea \nn &&\langle \bar
\kappa(\pi^p_x,\pi^a)-y_{\bar D},y_{\bar D}-y\rangle \\
\nn &=& \langle \bar \kappa(\pi^p_x,\pi^a)-x_{\bar D} + (x_{\bar
D}-y_{\bar D}) ,
y_{\bar D}-y\rangle \\
\nn &=& \langle \bar \kappa(\pi^p_x,\pi^a)-x_{\bar D}, x_{\bar D} -
x \rangle + \langle \bar \kappa(\pi^p_x,\pi^a)-x_{\bar
D},\left(y_{\bar D}-y-(x_{\bar D}-x)\right)\rangle
+ \\
\nn && \ \ \langle x_{\bar D}-y_{\bar D},
y_{\bar D}-y\rangle \\
\nn &>& \epsilon - |\langle \bar \kappa(\pi^p_x,\pi^a)-x_{\bar
D},\left(y_{\bar D}-y-(x_{\bar D}-x)\right)\rangle| - |\langle
x_{\bar D}-y_{\bar D},y_{\bar D}-y\rangle | \eea

Since $\sup_{\pi^a\in\Pi^a}\sup_{x_{\bar D}\in\bar D}\|\bar
\kappa(\pi^p_x,\pi^a)-x_{\bar D}\|<\infty$ and  $\sup_{y\in
K}\|y_{\bar D}-y\|<\infty$, it follows that there exists a finite
positive constant $c$ such that

\[
\langle \bar \kappa(\pi^p_x,\pi^a)-y_{\bar D},y_{\bar D}-y\rangle >
\epsilon - c(\|x_{\bar D}-y_{\bar D}\| + \|x-y\|).
\]

Since $\bar D$ is convex, the map $x\mapsto x_{\bar D}$ must be
continuous. It follows that there exists a $\rho(x)>0$ such that
$\langle \bar \kappa(\pi^p_x,\pi^a)-y_{\bar D},y_{\bar D}-y\rangle >
\epsilon/2$ whenever $\|x-y\|<\rho(x)$. Since this holds for any
$\pi^a$, we get
\[
 \inf_{y\in B(x,
\rho(x))}\inf_{\pi^a\in\Pi^a}\langle \bar
\kappa(\pi^p_x,\pi^a)-y_{\bar D},y_{\bar D}-y\rangle
> 0.
\]
\ep

For the rest of the paper we assume that $\rho(\cdot) : K\backslash
\bar D \ra \mathbb R_+$ is a function
satisfying~(\ref{satisfactory}). We now introduce the main objects
needed for our analysis.

\emph{The sets $K_n$, $\mathcal Q_n$ and $\mathcal Q$.}\ \  For
$n\in\mathbb N$, define compact sets $K_n$ as
\[
K_n := \left\{y\in K:\inf_{x\in D}\|y-x\|\in [1/(n+1), 1/n]\right\}.
\]
We can write
\[
K\backslash \bar D = \bigcup_{n\in\mathbb N} K_n.
\]
For $n\in\mathbb N$, the collection $\{B(x, \rho(x)/2): x\in K_n\}$
is an open cover for $K_n$. By compactness there exists a finite
subcover. Let $\mathcal Q_n$ be a finite subset of $K_n$ such that
\[
\bigcup_{q\in\mathcal Q_n}B(q, \rho(q)/2) \supset K_n.
\]
Let $\mathcal Q$ denote the union
\[
\mathcal Q := \bigcup_{n\in\mathbb N} \mathcal Q_n.
\]
The following result is immediate.
 \bprop

The collection $\mathcal Q$ is a countable collection.

\eprop

\emph{The map $Q(\cdot)$.}\ \  Since $\mathcal Q$ is countable, we
can assign an injective (one--one) map $I:\mathcal Q \lra \mathbb
N$. Using the map $I(\cdot)$ we define a map $Q:K\backslash \bar D
\ra\mathcal Q$ where, for $x\in K\backslash \bar D$, we define
\[
Q(x) := \underset{q}{\operatorname{argmin}}\{I(q): x \in B(q,
\rho(q)/2), q\in \mathcal Q\}.
\]

\emph{The re--scaled times and the interpolated trajectory.}\ \  Let
$t(0)=0$. For $n\in \mathbb N$, define the re--scaled times
\[
t(n)=\sum_{i=1}^{n}1/i.
\]
Let $\bar x(\cdot)$ be the trajectory obtained by linearly
interpolating between the iterates. Thus, for any $n \in \mathbb N$
and $t\in[t(n), t(n+1))$ define
\[
\bar x(t) := \frac{t(n+1)-t}{t(n+1)-t(n)}\cdot x_n +
\frac{t-t(n)}{t(n+1)-t(n)}\cdot x_{n+1}.
\]

\emph{The map $T(\cdot)$.}\ \  Define $v_{\max}:= \sup_{x\in
K}\sup_{(\theta,u^p,u^a)}\|\kappa(\theta,u^p,u^a)-x\|$. Since $x,
\kappa(\theta,u^p,u^a)\in K$ and $K$ is compact, it follows that
$v_{\max}<\infty$. Clearly, for times $u_1$ and $u_2$,
\[
\|\bar x(u_1) - \bar x(u_2)\|\leq v_{\max}|u_1-u_2|.
\]
Let $T:\mathcal Q\ra \mathbb R_+$ be a map such that for every $q
\in \mathcal Q$ the following holds:

\be \label{equationT} \frac{\rho(q)}{4 v_{\max}} < T(q) <
\frac{\rho(q)}{3v_{\max}}. \ee

\emph{Choice of strategy along $\mathcal S$.} \ \ We are now ready
to define how the player should choose his strategies over time. Let
$\pi^p_0$ be any arbitrary strategy. Let $\mathcal S:=(s_n)$ denote
the increasing subsequence of times when the player changes his
strategy. Start with $s_0=0$. Assume $s_n$ is known. We consider two
cases, $x_{s_n}\in K\backslash \bar D$ and $x_{s_n}\in \bar D$. If
$x_{s_n}\in K\backslash \bar D$ then set $q=Q(x_{s_n})$. Now choose
the strategy $\pi^p_q$ and set
\[
s_{n+1} = \underset{m}{\operatorname{argmin}}\left\{m:
\sum_{i=s_n}^{m-1}1/i > T(q)\right\}.
\]
If, however, $x_{s_n}\in \bar D$ then choose the strategy $\pi^p_0$
and set $s_{(n+1)}=s_n+1$.

\section{A two time scale result}
\label{tts}

This section develops an elementary two time scale result needed for
the proof of convergence. For the reader's convenience we break the
proof into a series of smaller units.

\bl \label{lemmarho} For every $x \in K \backslash \bar D$,

\[
B(x, \rho(x))\cap\bar D = \emptyset.
\]

\el

\bp

If $y\in B(x, \rho(x))\cap\bar D$, then $y_{\bar D}=y$ and so
\[
\inf_{\pi^a\in\Pi^a}\langle \bar \kappa(\pi^p_x,\pi^a)-y_{\bar
D},y_{\bar D}-y\rangle = 0.
\]
But this contradicts~(\ref{satisfactory}). Hence, $B(x,
\rho(x))\cap\bar D = \emptyset$.

\ep

\bl

\label{lemmafinite} For any compact set $L$ such that $L\bigcap \bar
D = \emptyset$, we have
\[
\left|\left\{q\in\mathcal Q: B(q, \rho(q)/2)\bigcap L \neq
\emptyset\right\}\right|<\infty.
\]
 \el

\bp

Since both $L$ and $\bar D$ are compact sets, it follows that
\[
\inf_{x\in\bar D, y\in L} \|x-y\| =: d(L) > 0.
\]

Consider any $q$ such that $q\in \mathcal Q_m$ and $m > \frac 3
{2d(L)}$. Since $\mathcal Q_m \subset K_m$, we have
\[
\inf_{x\in \bar D}\|q-x\| \leq \frac 1 m < \frac{2d(L)}3.
\]
Further, by Lemma~\ref{lemmarho}, $\rho(q)<1/m< 2d(L)/3$. It follows
that if $m > \frac 3 {2d(L)}$ and $q\in \mathcal Q_m$ then $B(q,
\rho(q)/2)\bigcap L=\emptyset$. The result follows.

\ep

\bl \label{lemmaintermediate}

Let $(s_{m(n)})$ be an increasing subsequence of $\mathcal S$. If
$\lim_{n\ra\infty} x_{s_{m(n)}} = x$ for some $x\in K\backslash \bar
D$, then along a further subsequence, denoted $(s_{m(n)})$ again,
there exists $q \in \mathcal Q$ such that $Q(x_{s_{m(n)}}) = q$ for
all $n \in \mathbb N$.

\el

\bp

Since $\lim_{n\ra\infty} x_{s_{m(n)}} = x \notin \bar D$, there
exists a compact set $L$ such that $L\cap\bar D=\emptyset$ and
$x_{s_{m(n)}}\in L$ for all sufficiently large $n$. By
Lemma~\ref{lemmafinite},
\[
\left|\left\{Q(x_{s_{m(n)}}): n\in\mathbb N\right\}\right|<\infty.
\]
Thus there exists $q \in \mathcal Q$ such that along a subsequence,
denoted $(s_{m(n)})$ again, we have $x_{s_{m(n)}}\in L$ and
$Q(x_{s_{m(n)}}) = q$ for all $n \in \mathbb N$.

\ep

\emph{The Mannor-Tsitsiklis bound.} \ \ We now introduce a set of
conditions, labeled $(\dag)$, which is needed for
Theorem~\ref{ldthm} and Corollary~\ref{corollary} below. To this
end, let $(s_{m(n)})$ be an arbitrary increasing subsequence of
$\mathcal S$. Let $T_l$ and $T_r$ be times such that $T_l < T_r$.
Let $(l_{m(n)})$ and $(r_{m(n)})$ be sequences such that
$s_{m(n)}\leq l_{m(n)} < r_{m(n)} \leq s_{m(n)+1}$, $n\in\mathbb N$.
Let $(\dag)$ denote the following four conditions: \bi
\item[$1^\dag$] $x_{s_{m(n)}}\lra x$ for some $x \in K\backslash \bar D$.
\item[$2^\dag$] $Q(x_{s_{m(n)}}) = q$ for some $q\in \mathcal Q$ and all $n\in\mathbb N$.
\item[$3^\dag$] $[T_l, T_r) \subset [0, T(q))$.
\item[$4^\dag$] $t(l_{m(n)})-t(s_{m(n)})\ra T_l$ and
$t(r_{m(n)})-t(s_{m(n)})\ra T_r$.
\ei

Assuming the conditions of (\dag) hold, for $l_{m(n)}\leq j <
r_{m(n)}$ consider the single step reward $\kappa\left(\theta_j,
u^p_j, u^a_j\right)$. At each of these time steps the player adopts
the strategy $\pi^p_q$ independently of the action chosen by the
adversary. For $\theta \in S$ and $u^a \in U^a$, let
$\kappa^\tau(\theta, u^a)$ be the reward at the $\tau^\text{th}$
occurrence of $(\theta, u^a)$ in the range
$l_{m(n)},\ldots,r_{m(n)}-1$, . The rewards $\kappa^\tau(\theta,
u^a), \tau = 1,2,\ldots,$ are independent, identically distributed
random variables with mean

\[
\mathbb E[\kappa^\tau(\theta, u^a)] = \sum_{u^p}\kappa(\theta, u^p,
u^a)\pi^p_q(\theta)(u^p).
\]
Further, since each $\kappa\left(\theta_j, u^p_j, u^a_j\right)$ is
chosen from a compact set, we get, for $z$ in any neighbourhood of
the origin,
\[
\mathbb E[\exp{(\langle z, \kappa^\tau(\theta, u^a)\rangle)}] <
\infty,
\]
where $\langle\cdot,\cdot\rangle$ is the inner product in $\mathbb
R^d$. Define the set $R(\pi^p_q):=\{\bar
\kappa(\pi^p_q,\pi^a):\pi^a\in\Pi^a\}$. For a vector $v$, define
$\|v - R(\pi^p_q)\|:=\inf_{\pi^a\in\Pi^a}\|v-\bar
\kappa(\pi^p_q,\pi^a)\|$. We can now invoke Theorem 6.2 of Mannor
and Tsitsiklis \cite{mannor}. For our setup and with our notation,
it reads as follows:

\bt \label{ldthm}

Assuming that the conditions of $(\dag)$ hold, there exists a
function $\lambda : (0, \infty) \rightarrow (0, \infty]$ and a
positive constant $c_0$, such that irrespective of the adversary
policy $\pi^a$, the following bound holds:
\[
\mathbb
P\left[\left\|\frac{\sum_{j=l_{m(n)}}^{r_{m(n)}-1}\kappa\left(\theta_j,
u^p_j, u^a_j\right)}{r_{m(n)} - l_{m(n)}} -
R(\pi^p_q)\right\|\geq\epsilon\right] \leq
c_0\exp{\left(-\lambda(\epsilon)(r_{m(n)} - l_{m(n)})\right)}.
\]
\et

For the next result, note that $t(r_{m(n)})-t(l_{m(n)}) =
\sum_{j=l_{m(n)}}^{r_{m(n)}-1}1/j$. Under the conditions of $(\dag)$
this implies that \be \label{equationrl} \lim_{n\rightarrow \infty}
\frac{r_{m(n)}}{l_{m(n)}} = \exp{(T_r - T_l)}. \ee

\bc \label{corollary} Assuming that the conditions of $(\dag)$ hold,
we have
\[
\lim_{n\rightarrow\infty}\left\|\frac{\sum_{j=l_{m(n)}}^{r_{m(n)}-1}\kappa\left(\theta_j,
u^p_j, u^a_j\right)}{r_{m(n)} - l_{m(n)}} - R(\pi^p_q)\right\| = 0
\text{ a.s.}
\]
\ec

\bp Since $\lim_{n\rightarrow \infty}(t(r_{m(n)})-t(l_{m(n)}))=T_r -
T_l
> 0$, it follows from~(\ref{equationrl}) that $r_{m(n)} - l_{m(n)} >
l_{m(n)}[\exp{(T_r - T_l)}-1]/2$ for $n$ sufficiently large. Since
$l_{m(n)}\geq n$, we get $r_{m(n)} - l_{m(n)} > n[\exp{(T_r -
T_l)}-1]/2$ for $n$ sufficiently large. Plugging this estimate in
Theorem~\ref{ldthm} and noting that the constant $\epsilon$ is
arbitrary, a standard application of the Borel--Cantelli argument
gives the result. \ep

\emph{The two time scale result.} \ \ With Corollary~\ref{corollary}
available for use, we are ready for our main two time scale result.
Thus, let $(s_{m(n)})$ be an arbitrary increasing subsequence of
$\mathcal S$. Assume that $\lim_{n\ra\infty}x_{s_{m(n)}}= x$ for
some $x \in K\backslash \bar D$. By Lemma~\ref{lemmaintermediate}
there exists a $q\in \mathcal Q$ such that along a subsequence,
denoted again by $(s_{m(n)})$, $Q(x_{s_{m(n)}}) = q$ for all
$n\in\mathbb N$. For $n\in\mathbb N$ and $t \geq 0$ define the
trajectories \be \label{eqtraj}
\bar y_{m(n)}(t) := \left\{ \ba{lcl} \bar x(t(s_{m(n)})+t) & \text{ if } & t \leq t(s_{m(n)+1}) - t(s_{m(n)}) \\
\bar x(t(s_{m(n)+1})) & \text{ if } & t > t(s_{m(n)+1}) -
t(s_{m(n)}) \ea \right. \ee By the Arzela-Ascoli theorem there
exists a continuous trajectory $\bar y(\cdot)$ such that along a
subsequence, denoted again by $(m(n))$, $\lim_{n\ra\infty}\bar
y_{m(n)}(\cdot)=\bar y(\cdot)$ in the topology of uniform
convergence over compacts.

Set $T = T(q)$. For $k \in \mathbb N$ and $J_k = \{0, 1, \ldots,
2^k-1\}$ consider the finite collection of intervals \be
\label{countable} \mathcal C_k:=\left\{\left[2^{-k}jT,
2^{-k}(j+1)T\right):j\in J_k \right\}.\ee For $j\in J_k$, define
\[
\kappa_{j,k} := \frac{\bar y(2^{-k}(j+1)T) - \exp{(-2^{-k}T)}\bar
y(2^{-k}jT)}{1-\exp{(-2^{-k}T)}}.
\]
Next, with $(s_{m(n)})$ denoting an arbitrary increasing subsequence
of $\mathcal S$, define $N(q,j,k)$  as the following set:
\[
N(q,j,k):=\{(x_n): \exists (s_{m(n)}) \text{ s.t. }
Q(x_{s_{m(n)}})=q\ \forall n\text{ and }
\kappa_{j,k}\notin\overline{R(\pi^p_q)}\}.
\]

\bprop

The set $N(q,j,k)$ is a null set, i.e., $\mathbb P[N(q,j,k)] = 0$.

\eprop

\bp

Fix any interval $\left[2^{-k}jT, 2^{-k}(j+1)T\right)$ in $\mathcal
C_k$. Let $(l_{m(n)})$ and $(r_{m(n)})$ be sequences with
$s_{m(n)}\leq l_{m(n)} < r_{m(n)} \leq s_{m(n)+1}$, $n \in \mathbb
N$ such that $t(l_{m(n)})-t(s_{m(n)})\ra 2^{-k}jT$ and
$t(r_{m(n)})-t(s_{m(n)})\ra 2^{-k}(j+1)T$. We have
\[
\frac{\bar y(2^{-k}(j+1)T) - \bar y(2^{-k}jT)}{2^{-k}T} =
\lim_{n\ra\infty}\frac{x_{r_{m(n)}}-x_{l_{m(n)}}}{t(r_{m(n)})-t(l_{m(n)})}.
\]

In terms of $l_{m(n)}$ and $r_{m(n)}$, the equation for average
reward can be written as
\[
 x_{r_{m(n)}} = x_{l_{m(n)}} +
\frac{r_{m(n)}-l_{m(n)}}{r_{m(n)}}\left[\frac{\sum_{j=l_{m(n)}}^{r_{m(n)}-1}\kappa\left(\theta_j,
u^p_j, u^a_j\right)}{r_{m(n)}-l_{m(n)}} - x_{l_{m(n)}}\right].
\]
Rearranging, we get
\[
\frac{x_{r_{m(n)}} -
(l_{m(n)}/r_{m(n)})x_{l_{m(n)}}}{1-(l_{m(n)}/r_{m(n)})} =
\frac{\sum_{j=l_{m(n)}}^{r_{m(n)}-1}\kappa\left(\theta_j, u^p_j,
u^a_j\right)}{r_{m(n)}-l_{m(n)}}.
\]

Since $x_{r_{m(n)}}\ra \bar y(2^{-k}(j+1)T)$ and $x_{l_{m(n)}}\ra
\bar y(2^{-k}jT)$, it follows from~(\ref{equationrl}) that
\[
\lim_{n\ra\infty} \frac{x_{r_{m(n)}} -
(l_{m(n)}/r_{m(n)})x_{l_{m(n)}}}{1-(l_{m(n)}/r_{m(n)})} =
\kappa_{j,k},
\]
and consequently
\[
\lim_{n\ra\infty}\frac{\sum_{j=l_{m(n)}}^{r_{m(n)}-1}\kappa\left(\theta_j,
u^p_j, u^a_j\right)}{r_{m(n)}-l_{m(n)}} = \kappa_{j,k}.
\]
Hence, by~Corollary~\ref{corollary} it must be the case that
\[
\kappa_{j,k} \in\overline{R(\pi^p_q)} \text{ a.s. }
\]

\ep

Define $\mathcal C:= \bigcup_k \mathcal C_k$. The next fact is
crucial to our analysis.

\bprop The collection $\mathcal C$ is a countable collection of
intervals.\eprop

Define $N$ to be the following set:

\[
N:=\bigcup_{q \in \mathcal Q}\bigcup_{k \in \mathbb N}\bigcup_{j \in
J_k} N(q, j, k).
\]

\bprop The event $N$ is a null set, i.e., $\mathbb P[N] = 0$. \eprop

\bp \label{null} Both $\mathcal Q$ and $\mathcal C$ are countable
collections. The result now follows from the fact that the union of
countably many exceptional null sets is again a null set.

\ep

By virtue of Proposition~\ref{null}, to show almost sure convergence
of sequences $(x_n)$ to $\bar D$ it suffices to restrict attention
to sequences outside $N$. Consequently, in what follows we shall
work exclusively with sequences $(x_n)$ outside the exceptional null
set $N$.

\bt Let $(x_n)$ be any sequence outside the exceptional null set
$N$. For $(s_{m(n)})$ an increasing subsequence of $\mathcal S$,
assume that  $\lim_{n\ra\infty} x_{s_{m(n)}} = x$ for some $x\in
K\backslash \bar D$. Assume further that for some $q\in\mathcal Q$,
$Q(x_{s_{m(n)}}) = q$ for all $n\in\mathbb N$. Let $T=T(q)$. Let
$\bar y(\cdot)$ be a limiting trajectory of the trajectories $\bar
y_{m(n)}(\cdot)$ given by~(\ref{eqtraj}). Then, for $t\in [0,T]$,
$\bar y(t)$ can be written as \be \label{twoeq} \bar y(t) = \bar
y(0) + \int_0^t v(s) \mathrm ds, \ee where $v(\cdot)$ is a Borel
measurable function defined on $[0,T]$. Further, for Lebesgue almost
all $t$ in $[0,T]$, the following holds: \be \label{derivative}
v(t)+\bar y(t)\in\overline{R(\pi^p_q)}. \ee
 \et

\emph{Remark.}\ \ We point out that~(\ref{twoeq}) is a standard
result in two time scale theory. Moreover, using Lebesque's theorem
we could also show~(\ref{derivative}) to hold almost surely for any
(but not all) $t\in[0,T]$. The problem arises from the fact that the
set $[0,T]$ is an uncountable set and when we do a union of null
sets, one for each $t\in[0,T]$, the union need not be a null set. We
solve this problem by treating the interval $[0,T)$ as a probability
space and giving the trajectory $\bar y(t)$ a martingale structure.
This also provides an independent and elementary proof of two time
scale structure.

\bp

Define $\mathcal G_k:= \sigma\left(\mathcal C_k\right)$, the
$\sigma$-algebra on $[0, T)$ generated by $\mathcal C_k$. Let
$\mathcal G := \bigvee_k \mathcal G_k$. For $\lambda$ the Lebesgue
measure, define the scaled probability measure $\mu$ on $[0, T)$
given by $\mathrm d\mu/\mathrm d\lambda = 1/ T$. This acts as a
probability measure for the probability space $([0, T), \mu,
\mathcal G)$. For $t\in[0,T)$ and $k\in\mathbb N$ define the `floor'
$f_k(t) := 2^{-k}\lfloor 2^kt/T\rfloor T$. Thus, for any
$t\in[0,T)$, we have $t\in [f_k(t), f_k(t)+2^{-k}T)$. Define
$M_k(t)$ as:
\[
M_k(t):= \frac{\bar y(f_k(t)+2^{-k}T) - \bar y(f_k(t))}{2^{-k}T}.
\]

Note that $M_k(\cdot)$ is $\mathcal G_k$--measurable. Further, for
$t\in[0, T)$ we have
\[
\mathbb E^\mu[M_{k+1}(t)| \mathcal G_k] = M_k(t) \ \mu\text{--almost
surely}.
\]
In other words, the sequence $(M_k(\cdot))_{k\in\mathbb N}$ forms a
bounded martingale in the filtered probability space $([0, T), \mu,
\mathcal G, \mathcal G_k)$. It follows that $\mu$--almost surely the
limit $v(t):= \lim_{k\ra\infty}M_k(t)$ exists. The limit $v(\cdot)$
is, clearly, a measurable function. Note that $[0, f_k(t))$ is a
$\mathcal G_k$--measurable subset of $[0, T)$. Letting $A:=[0,
f_k(t))$, it is immediate that $\int_A M_k(s)\mathrm ds = \int_A
v(s)\mathrm ds$. It follows that
\[
\bar y(f_k(t)) = \int_0^{f_k(t)}M_k(s)\mathrm ds =
\int_0^{f_k(t)}v(s)\mathrm ds.
\]
Letting $k\ra\infty$ gives us:
\[
\bar y(t) = \bar y(0) + \int_0^t v(s) \mathrm ds.
\]

Let $t\in [0, T)$. Set  $j = j(k) = \lfloor 2^kt/T\rfloor$. Note
that as $t$ ranges over $[0,T)$, the pair $(j(k),k)$ still take
values in a countable set. From the definitions of $M_k(t)$ and
$\kappa_{j,k}$ it follows that
\[
v(t)=\lim_{k\ra\infty}M_k(t) = \lim_{k\ra\infty}\kappa_{j(k),k} -
\bar y(t).
\]
Since $(x_n)$ is outside the exceptional null set $N$,
$\lim_{k\ra\infty}\kappa_{j(k),k}$ must necessarily lie in
$\overline{R(\pi^p_q)}$.

\ep

\bl \label{lemcrucial} Let $(x_n)$ be any sequence outside the
exceptional null set $N$. For $(s_{m(n)})$ an increasing subsequence
of $\mathcal S$, assume that  $\lim_{n\ra\infty} x_{s_{m(n)}} = x$
for some $x\in K\backslash \bar D$. Assume further that for some
$q\in\mathcal Q$, $Q(x_{s_{m(n)}}) = q$ for all $n\in\mathbb N$. Let
$T=T(q)$. Let $\bar y(\cdot)$ be a limiting trajectory of the
trajectories $\bar y_{m(n)}(\cdot)$ given by~(\ref{eqtraj}). Then
\[
\inf_{w \in \bar D}\|\bar y(t)-w\|\leq \inf_{w \in \bar D}\|\bar
y(0)-w\|\exp{(-t)} \text{ for all } t\in [0, T(q)].
\]
\el

\bp For $t\in [0, T)$ let $d(t):=\inf_{w \in \bar D}\|\bar
y(t)-w\|$. For any point $p$, let $d_p(t):=\|\bar y(t)-p\|$. Let
$\bar y_{\bar D}(t)$ be the point in $\bar D$ closest to $\bar
y(t)$. We have
\[
\dot{d}(t)\leq \dot{d}_p(t)|_{p = \bar{y}_{\bar{D}}(t)} =
\frac{(\bar{y}(t) - \bar{y}_{\bar{D}}(t))\cdot v(t)}{\|\bar{y}(t) -
\bar{y}_{\bar{D}}(t)\|} <-d(t),
\]
and the result follows.\ep

\section{Almost sure convergence}
\label{asc}

As before we present our proof as a series of short lemmas.

\bl

\label{lemmaquickchange} For $s_n\in\mathcal S$, if $x_{s_n}\in\bar
D$ then
\[
t(s_{n+1}) - t(s_n) = 1/(s_n+1) < 1/s_n,
\]
while if $x\in K\backslash \bar D$ then
\[
t(s_{n+1}) - t(s_n) < T(Q(x)) + 1/s_n.
\]
\el

\bl \label{lemmaTconvergence} For $x\in K\backslash \bar D$ we have
\[
T(Q(x)) \lra 0 \text{ as } x \lra \bar D.
\]
\el

\bp By definition, $x \in B(Q(x), \rho(Q(x))/2)$. It follows from
Lemma~\ref{lemmarho} that $B(Q(x), \rho(Q(x)))\bigcap \bar D =
\emptyset$. Consequently $\rho(Q(x)) \leq 2 \inf_{y\in\bar
D}\|x-y\|$. The result now follows from~(\ref{equationT}). \ep

\bl \label{lemmazero} Let $(s_{m(n)})$ be an increasing subsequence
of $\mathcal S$. If $\lim_{n\ra\infty} x_{s_{m(n)}} = x$ for some
$x\in\bar D$ then
\[
t(s_{m(n)+1}) - t(s_{m(n)}) \lra 0.
\]
\el

\bp By Lemma~\ref{lemmaquickchange}, if $x_{s_{m(n)}} \in \bar D$
then $t(s_{m(n)+1}) - t(s_{m(n)}) < 1/s_{m(n)} \leq 1/n$, while if
$x_{s_{m(n)}} \in K\backslash \bar D$ then $t(s_{m(n)+1}) -
t(s_{m(n)}) < T(Q(x_{s_{m(n)}}))+ 1/s_{m(n)}$.
By~Lemma~\ref{lemmaTconvergence}, $T(Q(x_{s_{m(n)}}))\ra 0$ as
$n\ra\infty$. Since $s_{m(n)}\ra \infty$ as $n\ra\infty$, the result
follows . \ep

\bl \label{lempair2} Let  $(s_{m(n)})$ be an increasing subsequence
of $\mathcal S$ such that $\lim_{n\ra\infty} x_{s_{m(n)}} = x$ and
$\lim_{n\ra\infty} x_{s_{m(n)+1}} = y$. If $y\in K\backslash\bar D$
then $x\in K\backslash\bar D$.

\el

\bp Assume $x\in\bar D$. Since $\|x_{s_{m(n)+1}} - x_{s_{m(n)}}\|
\leq v_{\max} (t(s_{m(n)+1})-t(s_{m(n)}))$, it follows from
Lemma~\ref{lemmazero}  that
\[
\lim_{n\ra\infty} \|x_{s_{m(n)+1}} - x_{s_{m(n)}}\| =0.
\]
This leads to a contradiction since  $y\in K\backslash \bar D$.
\ep

Recall that $\mathcal S= (s_n)$ is the increasing sequence of times
when the player changes his strategy.

\bl \label{llimpoint} Let $(x_n)$ be a sequence outside the
exceptional null set $N$. If $y$ is a limit point of the sequence
$(x_{s_n})$ then $y \in \bar D$. \el

\bp Assume to the contrary and let $y$ be a limit point of
$(x_{s_n})$ that is farthest from $\bar D$. Take an appropriate
subsequence such that $\lim_{n\ra\infty}x_{s_{m(n)}}=x$ and
$\lim_{n\ra\infty}x_{s_{m(n)+1}}=y$. By Lemma~\ref{lempair2} $x\in
K\backslash\bar D$. Further assume, by
Lemma~\ref{lemmaintermediate}, that the subsequence is such that
$Q(x_{s_{m(n)}})=q$ for some $q\in\mathcal Q$ and all $n\in\mathbb
N$. From our choice of $y$ it follows that
\[
\inf_{w\in\bar D}\|y-w\|\geq\inf_{w\in\bar D}\|x-w\|.
\]
But by Lemma~\ref{lemcrucial} we get
\[
\inf_{w\in\bar D}\|y-w\|\leq\inf_{w\in\bar D}\|x-w\|\exp{(-T(q))}.
\]
Since $T(q)>0$ this leads to a contradiction and the result follows.
 \ep

\bt \label{theoremmain} Let $(x_n)$ be a sequence outside the
exceptional null set $N$. If $x$ is a limit point of the sequence
$(x_n)$ then $x \in \bar D$. \et

\bp By taking suitable subsequences assume that
$\lim_{n\ra\infty}x_{u_{m(n)}}= x$ where $s_{m(n)} < u_{m(n)} \leq
s_{m(n)+1}$ for all $n\in\mathbb N$ with $(s_{m(n)})$ some
increasing subsequence of $\mathcal S$. Assume further that
$\lim_{n\ra\infty}x_{s_{m(n)}}= y$ for some $y$. By
Lemma~\ref{llimpoint} $y\in\bar D$. Since $\|x_{u_{m(n)}} -
x_{s_{m(n)}}\| \leq v_{\max} (t(s_{m(n)+1})-t(s_{m(n)}))$, it
follows from Lemma~\ref{lemmazero} that $\lim_{n\ra\infty}
\|x_{u_{m(n)}} - x_{s_{m(n)}}\| =0$. Thus $x=y$ and the result
follows.

\ep

\section{Conclusion}

\label{conclusion}

We have established the a.s.\ convergence of our scheme to the
desired limit set for finite state controlled Markov chains. In
conclusion we point out some future directions.

\emph{Extension to non--convex $D$.}\ \  For non-convex $D$ in
general, the existence of a `nearest point' in $\bar D$ from any
point outside $\bar D$ is guaranteed. A scheme along above lines can
be conceived wherein one uses piecewise constant policies that
ensure decrease of distance from $\bar D$ if such policies are known
to exist.

\emph{Countable state space.}\ \  Under suitable uniform stability
assumption or `near-monotonicity' condition on costs, variations of
the above scheme can be proposed for Blackwell approachability. This
will be pursued in a future work.

\emph{Computational issues.}\ \  The above scheme is an `ideal'
scheme in so far as it ignores actual computational aspects. A
practical implementation would raise further issues such as
recursive on-line computation of policies, learning, etc.

\emph{A combination scheme.} \ \ A  variation that seems promising
is to combine the approaches of this paper and Shimkin and
Shwartz~\cite{shimkin}, switching strategies when the currently
adopted strategy exhausts its allotted time, or when the chain
returns to a prescribed state, whichever occurs first. One expects
similar results, though the analysis will be messier.

\textbf{Acknowledgements.} The author would like to thank Prof.\ V.\
S.\ Borkar for introducing him to Blackwell approachability, for
pointing out reference~\cite{mannor} and for his careful reading of
an earlier draft and help with preparing this one.

\end{document}